\input amstex

\loadeufm
\loadmsbm
\loadeufm

\documentstyle{amsppt}
\input amstex
\catcode `\@=11
\def\logo@{}
\catcode `\@=12
\magnification \magstep1
\NoRunningHeads
\NoBlackBoxes
\TagsOnLeft

\def \={\ = \ }
\def \+{\ +\ }
\def \-{\ - \ }

\def \b|{\big |}

\def \g1{\Gamma_1}

\def\rarr#1#2{\smash{\mathop{\hbox to .5in{\rightarrowfill}}
 	 \limits^{\scriptstyle#1}_{\scriptstyle#2}}}

\def\larr#1#2{\smash{\mathop{\hbox to .5in{\leftarrowfill}}
	  \limits^{\scriptstyle#1}_{\scriptstyle#2}}}

\def\swarr#1#2 {\llap{$\scriptstyle #1$}  \swarrow
  	\vcenter to .5in{}\rlap{$\scriptstyle #2$}}

\topmatter
\title A report on Tarski's Decidability Problem:
\centerline{"Elementary Theory of Free Nonabelian  Groups"}
\centerline{by O. Kharlampovich and A. Myasnikov  }
\endtitle
\author
\centerline{
Z. Sela${}^{}$}
\endauthor
\abstract\nofrills{}
This paper  contains a list of  crucial mistakes and counterexamples to some of the main statements in the paper
"Elementary theory of free nonabelian groups" by O. Kharlampovich and A. Myasnikov, which
was published in Journal of Algebra in June 2006. 
\endabstract
\endtopmatter

\document

\baselineskip 12pt
O. Kharlampovich and A. Myasnikov announced a solution to Tarski's problems on the elementary theory of
the free groups in June 1998. Their work appears in a sequence of papers ending with the paper
 "Elementary theory of free nonabelian groups" that was
published in the Journal of Algebra in 2006 [KM4]. I had already written a report on this paper, reviewing 
the published version as well as  approximately 30 versions that preceded it, including serious mistakes
that appeared in essential points in all the versions.
In the current report I
single out only the (fatal) mistakes in the published paper, 
together with counterexamples to many of its main 
statements. 

The report starts with a short introduction that describes briefly the general approach to Tarski's problems that is 
presented in [Se1]-[Se7] and was adapted by the authors. As [Se1]-[Se7] mainly prove quantifier elimination,
and the arguments there are not effective, I further explain what needs to be proved, in order to construct
an effective procedure that will prove the decidability of the elementary theory of free
(or more generally all torsion-free hyperbolic) groups. 

The paper  continues with a short account of the main flaws/gaps in the paper by Kharlampovich and Myasnikov, that makes
it clear that no proof of any of Tarski's problems, in particular the decidability of the first order
theory of the free group, can be found in their paper. As these main flaws/gaps transfer directly to the recent paper
of the authors on the decidability of the theory of torsion-free hyperbolic groups [KM5], the report
clarifies that the decidability of the theories of both free and hyperbolic groups should be 
considered  as open problems.

The report  ends  with a detailed account of  mistakes and gaps in the paper under review (this list is essentially
contained in the previous report). The list is by no means 
complete, but it is sufficient in order to demonstrate and verify my arguments. 
Many of the indicated mistakes  are
followed by
counterexamples to the corresponding wrong claims in the paper.

\medskip
\centerline{\bf{1. The general approach towards Tarski's problems}}
\smallskip

The study of Tarski's problems starts with the analysis of sets of solutions to systems of equations over a free group.
The major breakthrough in the study of these systems of equations  is due to Makanin [Ma1] and Razborov [Ra]. Another important 
"classical" theorem
in the general direction of Tarski's problems is due to Merzlyakov from 1966 [Me]. Merzlyakov proved the
existence of "formal solutions" that essentially enable one to
analyze positive sentences, i.e., sentences that contain no inequalities.  
In their paper on what they call "implicit function theorem"  Kharlampovich and
Myasnikov present a generalization of Merzlyakov's theorem  to AE sentences (sentences with 2 quantifiers)
that are defined over general
varieties [KM3]. 

Merzlyakov's original theorem enabled Makanin to prove that the positive theory of the free groups is decidable [Ma2]. 
However, the
generalization of Merzlyakov's theorem to general (not only positive) sentences and formulas, enables one
to prove the validity of an AE sentence for a generic point of a variety and not for all the points in it.
The major problem and the main difficulty in tackling Tarski's problems is the ability to find a procedure that 
produces finitely many formulas that prove the validity of a sentence for all the points in a variety 
and not just for generic
points in it. This is the goal of most of our work on these problems [Se3]-[Se7], and what Kharlampovich and
Myasnikov are aiming to do in their paper under review [KM4].

As [KM4] uses exactly the same strategy, objects, constructions and procedures as [Se3]-[Se7], we briefly review the content of
[Se3]-[Se7], and use it to indicate what parts are missing or done wrongly in [KM4].
In [Se3] we studied the structure of exceptional solutions of a parametric 
system of equations (see definitions 10.5 in [Se1] and 1.5 in [Se3] for these exceptional solutions). 
We proved the existence of a global bound (independent of the 
specialization of the defining parameters) on the number of exceptional (rigid) solutions of a rigid limit
group ([Se3],2.5), and a global bound on the number of exceptional (strictly solid) families of solutions of a solid
limit group ([Se3],2.9). These global bounds are essential in our approach, and in particular are used to prove that
certain basic definable sets, that are later proved to be the "building blocks" of all the definable sets over a free
group, are  
in the Boolean algebra of $AE$ sets (section 3 in [Se3]).

In [Se4] general $AE$ sentences (sentences with two quantifiers) are analyzed.
Given a truth sentence of the form:  
$$\forall y \ \exists x \ \Sigma(x,y,a)=1 \, \wedge \, \Psi(x,y,a) \neq 1$$
we presented an iterative procedure,  that produces
a sequence of varieties and formal solutions defined over them, that together prove the validity 
of the given sentence. The procedure uses a trial and error approach. It starts with a formal solution 
(a formula) that proves the validity of the given sentence in a generic point of the affine set associated with
the corresponding universal ($y$) variables. If we substitute the given formal solution into the
system of equations $\Sigma$ they hold over the whole affine set. However, the inequalities $\Psi$
fail on some proper subvariety $V$. Hence, in the second step of the iterative procedure, we apply theorem
1.18 in [Se3] (what Kharlampovich and Myasnikov call "implicit function theorem"), 
and get formal solutions (formulas) that prove the validity of the given sentence in generic points of
certain Diophantine sets
(closures of completions) of resolutions in the Makanin-Razborov diagram of the variety $V$. Again, if
we substitute
the given formal solutions into the inequalities $\Psi$, 
they fail to hold on some subsets of the Diophantine sets that are associated 
the variety $V$. We continue iteratively by constructing
formal solutions (formulas) over a decreasing sequence of  Diophantine sets, that are defined using systems
of equations in increasing sets of variables.

Although in each step in the iterative procedure, one is forced to increase the set of variables, and even though
a decreasing sequence of Diophantine sets over a free group does not have to terminate, the procedure that is presented in [Se4] 
is constructed in such a way that a certain complexity of the varieties that are constructed along the procedure
strictly decreases, and forces the iterative procedure to terminate after finitely many steps.
It should be noted that even at this stage we don't know any conceptual reason
that explains the ability to find such a terminating procedure, except for its existence.

The outcome of the terminating iterative procedure for the analysis of an AE set (that is presented in [Se4])
is a collection of varieties, together
with a collection of formal solutions that are defined over them. The varieties are determined
by the original universal variables $y$, and extra (auxiliary) variables. 
The collection of varieties gives
a partition of the initial domain of the universal variables $y$, which is a power of
the original free group of coefficients, into sets which are in the
Boolean algebra of universal sets, so that on each such set the sentence can be validated
using a finite family of  formal solutions. Hence, the outcome of the iterative procedure can
be viewed as a $stratification$ $theorem$ that generalizes Merzlyakov's theorem from positive
sentences to general $AE$ ones.

\smallskip
In the two papers on quantifier elimination we apply the tools and techniques that are
presented in [Se1]-[Se4],
to prove quantifier elimination in the elementary theory of a free group. In 
order to prove quantifier elimination we show that the Boolean algebra of $AE$ sets is
invariant under projections. The projection of a set that is in the Boolean algebra of $AE$
sets, is naturally an $EAE$ set, hence, to show that the Boolean algebra of $AE$ sets is 
invariant under projections, we need to show that a general $EAE$ set is in the Boolean 
algebra of $AE$ sets ([Se6],1.41).  

To prove that an $EAE$ set is in the Boolean algebra of $AE$ sets we use a couple of
terminating iterative procedures that are based on the procedure for validation of an
$AE$ sentence presented in the fourth paper. 
Given an $EAE$ set, 
the first 
(terminating) iterative procedure is devoted to uniformization of proofs. i.e., it produces
finitely many (graded) families of formal solutions together with (graded) varieties on which these formal
solutions are defined, so that  each $AE$ sentence associated with a specialization of the defining
parameters and a specialization of the first existential variables, which is a truth sentence, can be proved using
part of the constructed families of formal solutions, in a similar way to our  validation of a (single)
 $AE$ sentence
obtained in the fourth paper.

Each step of this
procedure ([Se5],2) is divided into two parts. In the first part we collect all the formal solutions defined
over the (finitely many, graded) varieties that collect the set of those values of the universal variables,
for which the corresponding $AE$ sentence is yet to be proved. The second part uses the constructed formal 
solutions to get a proof for a subset of the relevant values of the universal variables, and collect
those values for which the proof is yet incomplete. We call the outcome of this
procedure, i.e., the families of formal solutions and the varieties on which they are defined, 
the $tree$ $of$ $stratified$ $sets$. Both its construction and its termination are 
uniformizations of the procedure for validation of a single $AE$ sentence that is presented in ([Se4],4).

The procedure for uniformization of proofs constructs the tree of stratified sets, that leaves us with
finitely many forms of proof, i.e., possible (finite) subsets of the families of formal solutions encoded
by this tree,   for all  the truth $AE$ sentences associated with the set
 $EAE$. we call each such form of proof a 
$proof$ $system$ 
([Se5],1.20).

To analyze an $EAE$ set, we start with the Zariski closures of all the $valid$ $proof$ $statements$  associated
with each of the (finitely many) proof systems ([Se5],1.23).
The second terminating iterative procedure that we call the $sieve$ procedure, presented in  [Se6], 
starts with each of these 
Zariski closures and constructs a (finite) sequence of bundles of (virtual)
proof statements  that are supposed to "testify" that a given specialization 
of the defining parameters is in the set $EAE$. This finite sequence of bundles reduces the question of the 
existence of a possible $witness$ (i.e., a value of the first existential variables)
with a valid proof statement ([Se5],1.19)
for any given specialization of the defining parameters, to the structure of the bases
of these bundles of proof statements. Since by section 3 of [Se3] it is possible to stratify the base
of such bundle, and the existence of a witness for a given specialization of the defining parameters 
depends only on the stratum (and not on the specific specialization), the set $EAE$ is the union of finitely
many strata in the stratifications of the constructed bundles.
Since every stratum in the stratification is in the
Boolean algebra of $AE$ sets ([Se3],3),
we are finally able to conclude that the original $EAE$ set is in the Boolean algebra of
$AE$ sets.

We should note that as in the procedure for validation of a sentence, presented in [Se4],
we still do not know a conceptual reason for obtaining quantifier elimination, and
for the ability to construct a terminating procedure like our sieve procedure, apart from its existence. Indeed,
the construction of the sieve procedure and its termination are technically the heaviest part of our work,
and require techniques and methods to handle Diophantine sets.

The quantifier elimination is the key for the analysis of the first order theory of a free group, and in particular
the key for solving Tarski's problems. The uniform quantifier elimination, i.e., a quantifier elimination that
does not depend on the rank of the (free) coefficient group, implies that all the non-abelian free groups are
elementarily equivalent - a solution to the first Tarski's problem. A possible approach to prove that the
first order theory of a free (or a hyperbolic) group is decidable (the second Tarski's problem), 
is proving that every step in our quantifier
elimination procedure can be made effective. This is precisely the approach that Kharlampovich and Myasnikov 
tried to apply in their attempt to prove the decidability of the theory. The quantifier elimination is also necessary and 
a key to later results, like proving the stability of the theories of free and hyperbolic groups, proving 
(geometric)
elimination of imaginaries, or analyzing the independence relation in the theory of a free group.

\medskip
\centerline{\bf{2. A general description of the crucial mistakes and gaps }}
\centerline{\bf{  in Kharlampovich and Myasnikov's paper}}
\smallskip
\smallskip

In the previous section we briefly sketched the approach to Tarski's problems as  presented
in [Se1]-[Se7]. In this section 2 we use this brief survey to indicate the main flaws and gaps in 
Kharlampovich and Myasnikov's
paper "Elementary theory of free nonabelian groups".  
A detailed account of the main flaws together with counterexamples to  claims in their paper 
appears in the next section. 
These flaws imply that none of Tarski's problems is proved
in this paper, and in particular, the decidability of the first order theory of both free and hyperbolic
groups, that is claimed in the paper and its more recent generalization to hyperbolic groups [KM5],
should be considered as open problems.

\roster
\item"{(1)}" The study of sets of solutions to systems of equations with parameters is basic and necessary
in any attempt to solve Tarski's problems or to obtain quantifier elimination. A fundamental property of
these parametric systems is the existence of a universal bound on the number of exceptional families of
solutions for such systems. This is the main goal of [Se3].

Kharlampovich and Myasnikov state a similar theorem in their paper, but the proof, although similar to the
one in [Se3], contains several serious mistakes. In particular, for effective results, it is essential to get an
effective bound on the number of exceptional families of solutions of parametric systems of equations. As
the proof of the existence of a bound in KM's paper is wrong, needless to say that the claim regarding an
effective bound is meaningless.

\item"{(2)}" In the procedure for the analysis of an AE sentence [Se4], and in the quantifier elimination procedure
[Se5]-[Se6], one of the main objects that is used repeatedly is the $induced$ $resolution$ that is presented 
in section 3 of [Se4].
The induced resolution is a canonical resolution, or a tower, that is inherited by a subgroup of an ambient tower.
Its construction is involved and uses an iterative procedure that is guaranteed to terminate.

Kharlampovich and Myasnikov are trying to construct the induced resolution and effectively compute it. However,
their construction is wrong. It does not produce a tower, and it can not serve for the purposes it was built for,
i.e., as a tool in the iterative procedure for analyzing definable sets. Once again, the construction is wrong, so
there is no point in considering the effective claim.

\item"{(3)}" AE sentences are analyzed in [Se4], using an iterative procedure that is guaranteed to terminate.
This procedure is still far from analyzing general definable sets, but it is essential in our approach for the
analysis of such general sets. For obtaining effective results, it seems necessary to be able to
run such a procedure effectively.

Kharlampovich and Myasnikov are trying to present exactly the same procedure, and claim to run it
effectively. However, their presentation is wrong,
and contains several fatal mistakes. What they describe may not terminate 
after finitely many steps.
Again, as the description of the procedure itself is wrong, there is no point in considering the effective parts.  

\item"{(4)}" The tools that are used in analyzing AE sentences are not sufficient for analyzing general
definable sets. One of the main new tools that we use for obtaining quantifier elimination is the
$core$ $resolution$, that is presented in section 4 of [Se5]. The construction of the core resolution is
rather complex and uses several iterative procedures. Any attempt to use the core resolution in effective 
procedures requires an effective computation of it.

Kharlampovich and Myasnikov are trying to present a construction of the core resolution. What they present 
contains several crucial mistakes and few major gaps, and hence, does not have the essential properties
 that  core resolutions must have  for the purposes of
quantifier elimination later in their paper. Again, the whole construction is wrong, so there is no point in 
referring to the effective claims.

\item"{(5)}" The heaviest machine or procedure that appears in our papers on Tarski's problems, is the sieve 
procedure for quantifier elimination [Se6]. This procedure is rather technical and complex, and is guaranteed
to terminate after finitely many steps, a termination that finally implies uniform quantifier elimination for all
non-abelian free groups (and a similar one for torsion-free hyperbolic groups).

Kharlampovich and Myasnikov aim to present precisely the same procedure, and run it effectively.
However, already some of the main objects that are used in the procedure, like the induced and the core
resolutions, were not constructed properly, and hence don't have the required properties for the procedure.
In addition, the description of the procedure contains several fatal mistakes and few gaps, so that the
presented procedure cannot be guaranteed to terminate. Again, there is no point in considering any of the
effective claims.
\endroster

\medskip
\centerline{\bf{3. A detailed description of the fatal mistakes and gaps }}
\centerline{\bf{ in Kharlampovich and Myasnikov's paper}}
\smallskip
\smallskip

In the previous section we briefly described some of the major flaws in Kharlampovich and Myasnikov's paper.
In this section we present a detailed account of some of these fatal mistakes and gaps. The description of the mistakes
is often followed by counterexamples to the statements in the paper under review. 
We should note that everything that appears in this section is contained in our previous report that
also describes the similarities with our papers.

\roster
\item"{(1)}" Page 470 lines 14-17. Lemma 8 is wrong. As stated, it is completely false. Even modulo inner automorphisms,
i.e., in $Out(G)$, the lemma is wrong as  was observed by G. Levitt, who found this mistake in an analogous theorem of mine
for torsion-free hyperbolic groups (from 1991). Levitt presented a counterexample and stated and proved a correct 
statement in [Le]. However, the 
(false) statement of the lemma 
does not play an essential role in the sequel. 

\item"{(2)}" Page 482 line -1 to page 483 line 2: Lemma 18 is false (lines 8-10 on page 483 give a wrong argument). No similar
statement can be correct, however, as the authors borrowed all the study of equations with parameters (a fully residually
free group relative to a subgroup) from [Se1] and [Se3] no such lemma is needed in the sequel. 

Below is a counterexample to the statement of lemma 18 (we use the notation of the lemma). Let $K$ be the limit group:
$$K \, = \, <u_1,\ldots,u_4> \, *_{[u_1,u_2][u_3,u_4]=[a^3,b^3]} \, <a,b> \, *_
{[a^4,b^4]=
[v_1,v_2][v_3,v_4]} $$
$$<v_1,\ldots,v_4>$$
and $K_1$ be the limit group:
$$K_1 \, = \, <u_1,v_1> \, *_{[u_1,v_1]=[a^5,b^5]} \, <a,b> $$

Let $\varphi:K \to K_1$ be the epimorphism that maps the tuple $u_3,u_4,v_3,v_4,u_2,v_2$ to the tuple:
$a^3,b^3,a^4,b^4,1,1$. Let $Q_1$ be the QH vertex group $<u_1,v_1>$, which is the unique QH vertex group in
the abelian JSJ decomposition of $K_1$. Then no QH subgroup in $K$ is mapped into a finite index subgroup in
$Q_1$, a contradiction to the conclusion of lemma 18 on page 482.

\item"{(3)}" Page 483 lines 11-16. Lemma 19 is false. Below is a counterexample to its statement.
Let: 
$$K \, = \, <a_1,a_2> \, *_{[a_1,a_2]=[b_1,b_2]} \, <b_1,b_2> \,
*_{[b_1,b_2]=[c_1,c_2]} \, <c_1,c_2> $$
$$H \, = \, <a_1,a_2> \, *_{[a_1,a_2]=[b_1,b_2]} \, <b_1,b_2>$$ 

$K$ has no sufficient splitting modulo $H$. Let $\tau: K \to H$ be the epimorphism that maps the tuple $c_1,c_2$ to 
$b_1,b_2$. By the notions of the authors, $H$ is a reducing quotient of $K$.
 No homomorphism  from $H$ (which is the fundamental group of a surface of genus 2) extends to a sequence of homomorphisms
that discriminate $H$, a contradiction to the conclusion of lemma 19.

\item"{(4)}" Page 485 lines -5 to -3. "since the group $K$ is ... one may assume that ... is a monomorphism". 
 
It was true for the definitions and the statement of their theorem 11 in some of their previous versions. However, in the
current statement of theorem 11, the authors
 "adapted" our notion of $strictly$ $solid$ homomorphisms (their "algebraic solutions"), and for this current statement one
can not make this reduction. 

It is worth noting that although the reduction that the authors make is wrong, 
a counterexample for this reduction is impossible to find, as a "growing sequence" does not
exist by theorems 2.5 and 2.9 in [Se3] (but in proving these theorems in [Se3], no reduction of the kind the authors
make is used).
 
\item"{(5)}" Page 486 lines 4-5: "Generalized equations for which $K$ is not embedded into ... result some of the 
equations from $R$".
The authors are making the same mistake as in page 485 lines -5 to -3. 
They (the authors) continue to refer to the statement of the 
theorem that appears in
their previous versions, and not to the current statement that uses our $strictly$ $solid$ homomorphisms ("algebraic
solutions") and $maximal$ $flexible$ $quotients$ ("complete reducing system").  

Page 486 line -7: "In case (2) we obtain an equation from the family $R$". The same mistake as in page 485 lines -5 to -3.
Note that once again, as "growing sequences" do not exist by theorems 2.5 and 2.9 in [Se3], and the authors study
these "growing sequences", no counterexamples can be constructed - still the arguments they use are wrong.

Theorem 11 (without its effective claim) is precisely theorem 2.9 in [Se1]. Hence, even though the argument the 
authors give
has several mistakes, the statement of the theorem (with the effective claim excluded) is correct. 
However, the statement of theorem 11 further claims that the bound
on the number of exceptional solutions can be effectively found. Since the argument that is given in the paper contains 
mistakes, the effective claim still requires a correct proof. It should be noted that the existence of an effective 
procedure to find a bound on the number of exceptional solutions of a parametric equation (the effective claim
 in the statement of theorem
11) is crucial for proving
decidability in the sequel.

\item"{(6)}" Page 506 line 10 to page 507 line 14: In this section the authors describe a construction that they
call "Induced NTQ system". 
What the authors aim to describe is the construction of our 
$induced$ $resolution$ (definition 
3.5 in [Se4]), but the construction that the authors describe contains several crucial mistakes.

Note that the induced resolution plays a major role in our papers, and is crucial in obtaining quantifier elimination.
In the sequel, wherever the authors use the $induced$ $resolution$ ("Induced NTQ system"), they are actually
using the object
that is constructed in section 3 of [Se4], and not the object that they are describing in their section 7.12. However,
for decidability one needs to effectively compute the induced resolution (as it appears in [Se4]). Because of the 
crucial mistakes
in the construction that the authors describe, no such effective computation is available in the current paper.

Below is an example for which the construction that the authors describe does not give an NTQ group at all (clearly
in contradiction to what they aim and claim to construct).
Let $K=<a,b> \, *_{[a^3,b^3]=[u,v]} \, <u,v>$. $K$ is clearly an NTQ group, where the base level is the
free subgroup $<a,b>$, and in the second level a fundamental group of a punctured torus $Q=<u,v>$ is amalgamated
to the base $<a,b>$ along its boundary subgroup. A fundamental group of a
punctured torus contains a subgroup of index 3, that we denote
$Q_1$, which is the fundamental group of a punctured torus with 3 punctures. Hence $K$ has a subgroup,
that we denote $G$, that has the structure:
$$G \, = \, ((<a_1,b_1> \, _{[a_1^3,b_1^3]=c_1} \, Q_1) *
_{[a_2^3,b_2^3]=c_2} \, <a_2,b_2>) \, *$$ 
$$*_{[a_3^3,b_3^3]=c_3} \, <a_3,b_3>$$
where $c_1,c_2,c_3$ are the 3 boundary components of $Q_1$. The procedure that the authors describe will not change 
the group $G$. But $G$ does not have the structure of an NTQ group.

\item"{(7)}"   Page 506  lines -11 to -10:
 "Increasing $G_1$ by a finite number of suitable elements from abelian vertex groups 
of $F_{R(S)}$ we join together ... from abelian vertex groups in $F_{R(S)}$". The sentence does not 
give any hint at what elements should be added. Could it be that the authors refer to what appears  in part (ii) in the
construction of the induced resolution in section 3 of 
[Se4] (starting at the 3rd paragraph after definition 3.2 in [Se4])?

\item"{(8)}" In all the construction of their "induced NTQ systems" the authors say nothing about QH vertex groups. 
In other words, they completely 
"skipped" part (iii) in our construction of the induced resolution (they have a previous section on "induced QH vertex groups" (lemma 7 on page 466)
that is not mentioned here. In any case, this notion and lemma are not a replacement of part (iii) in section 3 of [Se4]). 
Without this part (or any other part), the construction
that the authors describe doesn't  make any sense, and is therefore, completely wrong 
(once again, in all their referrals to "induced NTQ systems", 
the authors have to use an object that
is constructed according to the complete construction that appears in section 3 of [Se4], and not the vague, partial 
and mistaken
object that they are describing
in this section).

\item"{(9)}" On
page 507, lines 4 to 8, the authors aim to explain why the iterative procedure terminates after finitely
many steps (what is proved in proposition 3.3 in [Se4]). The 3 possible decreases they count are actually true. 
However, according to the proof of proposition 3.3 in [Se4], these decreases  happen on a highest level, and
on all the levels above that level the associated decompositions have the same structure as in the step before, 
and below that level there is no control on what happens with the associated decompositions (their complexity may increase).

The authors just indicate that in "some level" one of the possible decreases might occur. This does not imply termination
of the procedure. This mistake can be easily fixed according to the proof of proposition 3.3 in [Se4].

\item"{(10)}" Page 507 lines 9-10: "the image of the to $i$ levels of $F_{R_Q}$ on the level $j+1$ is the same as the 
image of $G$ on this level". This is completely false. It is true that the image is contained in the image of $G$ 
in this level (this is precisely the statement of lemma 3.6 in [Se4]). 

Below is a counterexample. Suppose that $K$ is the NTQ group: 
$$K \, = \, <u_1,u_2> \, *_{[u_1,u_2]=[a^3,b^3]} \, <a,b>$$ 
$K$ has a base group $<a,b>$ and in the second level the subgroup $<u_1,u_2>$, which is the fundamental group
of a punctured torus is amalgamated to $<a,b>$. 

Let $G=<u_1^2,u_2^2>$. Then the image of $G$ in the base level of $K$ is the subgroup $<a^6,b^6>$. However, $F_{R(Q)}$ is
a free group of rank 2, and its image in the base level is trivial, a contradiction to the claim of the authors that the
image of $F_{R(Q)}$ and $G$ in the base level are the same.
   
\item"{(11)}" The construction of the induced resolution in section 3 of [Se4] takes as an input a $well$-$structured$
resolution and a subgroup of its completion. The $well$-$structured$ structure is fundamental in the construction 
(see part (iii) in definition 3.5 in [Se4] that the authors completely skipped). 

A $well$-$structured$ structure is a slightly weaker assumption than a $well$-$separated$ structure (all these are our 
notions), that the authors "borrowed" in their sections 7.8 and 7.9 (pages 504-505) "first and second restrictions
on fundamental sequences". In section 7.10 "Induced NTQ systems" the "first and second requirements on fundamental
sequences" are never assumed or mentioned. Without this (geometric) structure, the whole construction 
of "induced NTQ systems" can not make
any sense  (and once again, when it is used in the sequel it's always according to the construction that appears in 
section 3 of [Se4], and (implicitly) under its assumptions and input). 

Therefore, no effective procedure that gets
as input a resolution (fundamental sequence) with a well-structured structure, and a subgroup of the corresponding
limit group, and computes effectively
the induced resolution is available in the current paper. Such an effective procedure is crucial for proving decidability.  

\item"{(12)}" Page 518 lines 20-22: "We may assume that $F_{R(U_0)}$ is freely indecomposable, otherwise we can 
effectively split into a free product of finitely many freely irreducible factors...and continue with each of the 
factors in the place of $F_{R(U_0)}$". This is a crucial mistake - a factorization followed by a reduction
that the authors suggest is absolutely forbidden. 
It demonstrates a basic misunderstanding of the procedure that is presented in section 4 of [Se4], and leads to
technical mistakes in the procedure itself.

Below is a counterexample. Consider the AE sentence:
$$ \forall y_1,y_2,y_3,y_4 \ \  [y_1,y_2] \neq 1  \, \vee  \, [y_3,y_4] \neq 1 \, \vee  \, y_1 \neq y_3 $$
$$ \vee \, y_2 \neq y_4 \, \vee [y_1,y_4]=1$$
The sentence is obviously true over a free group. However, after using the first two inequalities, the procedure
gives a freely decomposable subgroup: $<y_1,y_2,y_3,y_4> \, | \, [y_1,y_2], [y_3,y_4]>$ which is clearly freely
decomposable. Any attempt to apply any form of a generalized Merzlyakov's theorem (what the authors call "implicit
function theorem"), involves the two factors in the free decomposition. Hence, no reduction to separate factor(s) 
is possible,
a contradiction to the authors claim.

\item"{(13)}" Page 520 line -2 to page 521 line 1: "Let ... be the subset of homomorphisms... 
and satisfying the additional equation 
$U_1(X_1,\ldots,X_k)=1$". This is a crucial mistake that leads to further mistakes in the sequel. What one should consider
are only homomorphisms in $shortest$ $form$ as they are defined in definition 4.1 in [Se4]. 

\item"{(14)}" Page 521 lines 5-6: "modulo the images ... of the factors in the free 
decomposition of $H_1=<X_2,\ldots,X_m>$". This 
is what is done in the first step in section 4 of [Se4], but it it is true only if the authors would have considered
only $shortest$ $form$ homomorphisms (definition 4.1 in [Se4]). page 521 lines 21-22: "canonical sequences for $H_1$ 
modulo the factors in the free decomposition of the subgroup $<X_3,\ldots,X_m>$" - again the same mistake (true
only for $shortest$ $form$ homomorphisms, but these are not the homomorphisms that the authors consider).

\item"{(15)}" Page 521 lines -11 to -9: Lemma 23  is supposed to be identical to proposition 4.3 in [Se4],
which is a key observation for the construction and the termination 
of the iterative procedure for the analysis of AE sentences.  
However, because
the authors haven't used $shortest$ $form$ homomorphisms in the construction of their fundamental sequences, both
the formulation of lemma 23 and the  argument that is used for its proof (that imitates the proof of 
proposition 4.3 in [Se4], but the assumptions there are slightly different) are wrong.

In the formulation of the lemma, it is not true that "it is possible to replace $H_{(p)}$ by a finite number of
proper quotients of it without losing values of initial variables of $U=1$" (page 521 lines -10 to -9). 

The argument that is used to prove the lemma is very confused. The authors are trying to change homomorphisms
to be in $shortest$ $form$ (definition 4.1 in [Se4]), but it is somewhat late at this stage (they should have 
considered only such homomorphisms to start with). In particular, at this stage,
taking "minimal solutions with respect to 
$A_{D_t}$" (line -3 page 521) changes the values of "variables of $U=1$" and this is not allowed (i.e., by doing that 
one loses "values of initial variables of $U=1$" (line -2) that the procedure must handle).

The following is a counterexample to lemma 23. Let $L$ be the limit group:
$$L \, = \, <a_1,a_2> \, *_{[a_1,a_2]=[b_1,b_2]} \, <b_1,b_2> \, *_{[b_1,b_2]=[c_1,c_2]} \, <c_1,c_2>$$
This cyclic splitting is clearly the JSJ decomposition of $L$. $L$ admits the resolution ("fundamental sequence"
in the terminology of the authors): $L \to \pi_1(S_2)$, where $S_2$ is the fundamental group of an
orientable  surface of genus 2. Hence, $L$ naturally embeds into the $NTQ$ group $N$, that is the completion
of the given resolution, and is obtained from an
amalgamated product of $\pi_1(S_2)$ and $L$, along the edge groups in $L$ and a separating s.c.c.\ on $\pi_1(S_2)$:
$$N \, = \, L \, *_{[b_1,b_2]=[d_1,d_2]} \, <d_1,d_2> \, *_{[d_1,d_2]=[e_1,e_2]} \, <e_1,e_2>$$

Now, we impose a new relation on the NTQ group $N$, that the authors denote $U_1$. The new relation $U_1$ that we impose
is given by the relations: $$d_1(a_1)^{-1},d_2(a_2)^{-1},e_1(b_1)^{-1},e_2(b_2)^{-1}.$$ 
The group that is obtained from the NTQ limit group $N$ by adding these relations has several (maximal) limit quotients.
One of these maximal limit quotients is $L$ itself, where the map from the NTQ group $N$ onto $L$ is the natural
retraction. In our terminology this retraction $L$ is $solid$ with respect to the group that is generated
by $\pi_1(S_2)$.  This contradicts the conclusion of lemma 23 that claims that
in every such limit quotient, the original group $L$
does not embed.

\item"{(16)}" Page 522 lines 5-7" "block-NTQ group ... generated by the top p levels of the NTQ group corresponding
to the fundamental group c and the group ... corresponding to some branch of the tree $T_{CE}(G_{(p)})$". This is
supposed to be the construction of the $anvil$, as it is presented in part (2) of the first step of the procedure
in section 4 of [Se4], and in definition 4.5 in [Se4]. 

Note that what the authors write is mistaken, as in 
general one can not take an amalgamated product as they suggest, because the group $G_{(p)}$ 
that is associated with the terminal
level of their NTQ group may not be embedded in the group that is associated with the 
top level
of the branch of the tree $T_{CE}(G_{(p)})$. This technical difficulty is treated in part (2) of [Se4] but somehow
the authors missed it...  

Below is a simple counterexample. Let $L$ be the limit group that we used in comment (15):
$$L \, = \, <a_1,a_2> \, *_{[a_1,a_2]=[b_1,b_2]} \, <b_1,b_2> \, *_{[b_1,b_2]=[c_1,c_2]} \, <c_1,c_2>$$
With the terminology of [Se2], the strict Makanin-Razborov diagram of $L$ contains 3 resolutions that starts with an
epimorphism: $L \to \pi_1(S_2)$, where $S_2$ is an orientable surface of genus 2, and an additional resolution that
starts with a proper quotient of $L$, that is obtained from $L$ by killing the element $[b_1,b_2]$, and is
isomorphic to $Z^2*Z^2*Z^2$ (which is already a tower or a NTQ group in the terminology of the authors). 

In constructing the $anvil$, one should attach this last  group $Z^2*Z^2*Z^2$ to the previously constructed group, in
which the limit group $L$ is embedded. This can not be done by amalgamation (as the authors claim), 
as $Z^2*Z^2*Z^2$ is a proper quotient
of $L$ and is not isomorphic to it.
  
\item"{(17)}" Page 522 lines 8-10: "One can extract from c ... induced by the fundamental sequence c. Denote this 
$extracted$ $fundamental$ $sequence$ by $c_2$". What the authors are really refering to here is their
"Induced NTQ systems" (page 506 section 7.12).
However, whatever is constructed in section 7.12 can not work here (as what they constructed in section 7.12
is not even an NTQ group,
and no generalized form of Merzlyakov theorem (i.e., what the authors call implicit function theorem) 
can be applied to it). 

See the counterexample in part (6) - it explains why the authors' "Induced NTQ system" can not serve
for the purposes they need it here - the generalized form of Merzlyakov's theorem, what the authors call "implicit
function theorem", simply doesn't apply to it as in general it is not an NTQ group.

\item"{(18)}" Page 
522 lines 23 to 25: "Consider the set of those homomorphisms... and satisfy some additional equation $U_2=1$".
Once again, the same mistake as in the first step - only $shortest$ $form$ homomorphisms (definition 4.1 in [Se4])
should be considered. Omitting the shortest form requirement will cause difficulties in the sequel (as in the first
step).

\item"{(19)}" Page 
523 lines 1 to 5: "Case 2" - the same mistakes
that the authors made in their "first step" recurs here.

See the counterexample to lemma 23 in part (15). As the conclusion of lemma 23 is false, and Case 2 of the authors 
is based on the validity of this conclusion, Case 2 of the authors simply can
not be executed.

\item"{(20)}" Page 523 lines 6 to 11: "Case 3".
According to the authors, the case in which $G$ is mapped
to a proper quotient, and the case in which $N^1_0$ is mapped to a proper quotient are treated in the same way. 
This is
a crucial mistake. First, in both cases the authors use the conclusion of lemma 23, which is false (see part (15)).
Second, in part (4) of the second step in [Se4] these two cases are treated differently. The authors make the same 
mistake in the general case, and using this mistake it's possible to construct a non-terminating procedure (see comment
(21) below).

\item"{(21)}" Page 523 
line -11 to -6: "Case 2"  
The crucial mistake that the authors made in their construction of the "block NTQ group"
(i.e., the $anvil$) in Case 3 of the second step, recurs here. In lines -7 to -6
the authors construct an NTQ group "with the top part being $c^{(n)}$ above level $p$ and the bottom part $f_i$. 
This is wrong, and the argument that the authors use for the termination of their procedure (the proof of theorem
36 on page 524), which is the argument that is used in proving theorem 4.12 in [Se4], can not work for the
construction that the authors describe (it does work for the construction that appears in [Se4]).

In fact,  by slightly modifying
the counterexample of a non-terminating
procedure that we presented in 2001 [Se8] (as a counterexample to Theorem 5 in the original version of
the paper under review), it is not difficult to construct a concrete iterative
procedure of the type that the authors describe, that doesn't terminate
after finitely many steps. i.e., the whole analysis of AE sentences in the current paper is  wrong.

\item"{(22)}" Page 523 lines 
-16 to line -12. What is written there is completely wrong. 
If cases 2 and 3 do not apply at any level, 
one needs to look at the structure of the induced resolution.
By proposition 4.11 in [Se4], the complexity of the induced resolution is guaranteed to drop in this case.

The authors messed up all that. First, the authors do not consider the possibility
of a change in the induced resolution (what they call "induced NTQ system" in section 7.12 on page 506).
Second, what they do suggest is to replace the original group (that they denote as
$G$) by some sort of an envelope,
which is the image of our $developing$ $resolution$ that was constructed in the previous step. A change (an increase)
of the original  group by some sort of "envelope"  is what we do (under different assumptions,
and after checking the structure of the induced resolution!) in
the procedure that analyzes formulas with more than 2 quantifiers (i.e., the $sculpted$ $resolution$),
but not in this procedure that analyzes sentences with only 2 quantifiers. 
Here, no such "envelope"  is used in our procedure (there is no need to consider such), and of course,
there is no referral to anything like that later in proving the termination of the procedure (the whole proof
of theorem 36 on page 524 breaks down because of what the authors do in their "Case 1"). So the authors are 
"proving" (in their theorem 36)
the termination
of our procedure in section 4 of [Se4], but the procedure that they construct is completely different (and fatally wrong).
In other words, there is no reason why the procedure that the authors describe will actually terminate.

As in comment (21) it is not difficult to use the mistake in Case 1 of the authors (the replacement of the group $G$
by the group $G'$), to construct a procedure of the type that is described by the authors that does not terminate
after finitely many steps.
 
\item"{(23)}" Page 524 line 3. 
As we have already mentioned, the conclusion of theorem 36 is false (see comments (21) and (22)),
and the proof of theorem 36 
applies to the procedure in section 4 in [Se4], not to the (actual details of the)
procedure that the authors describe in their general step on page 523. The procedure that the authors describe 
does not terminate in general. It is not difficult to modify the non-terminating procedure in [Se8], to 
construct a counterexample to theorem 36 of the authors.

\item"{(24)}" Page  
525 lines 1 to 5:  The conclusion of lemma 24 is correct for the top level of the original fundamental sequence. 
However, lemma
23 is false (see comment (15)), so it makes it impossible to apply the constructions that appear in the
authors' procedure and apply lemma 24 in the next levels.

\item"{(25)}" The argument that is supposed to prove theorem 36 and appears on page 525 line -12 to page 526 line 3, that
actually briefly sketches the proof of theorem 4.12 in [Se4], applies to the procedure that is presented in section 4 
of [Se4], but not to the procedure for the construction of the tree $T_{AE}(G)$ that the authors presented
in pages 520-523. This is mainly because of the serious mistakes in Cases 1 and 2 of the authors' general step (see comments
21-23), and because the authors do not use homomorphisms in $shortest$ $form$ (definition 4.1 in [Se4]).

\item"{(26)}" Page 528 lines 19-21: "Consider a finite family of terminal groups of fundamental sequences of $P$ 
modulo factors in the free decomposition of  $F_{R(U)}$". This is a  fatal mistake. 
The diagram that the authors consider
is with respect to the freely indecomposable non-cyclic 
factors of $F_{R(U)}$. However, in a resolution (fundamental sequence in their
terminology) of such a diagram, the values of the free factor of $F_{R(U)}$ are being changed, and the values of the
other factors are modified by conjugation. Hence, not every value of
$F_{R(U)}$ that can be extended to a value of $F_{R(P)}$ can be extended to a value of one of the terminal groups
in the diagram that the authors are constructing. For example, it may be that $N_1$ is embedded in $F_{R(P)}$ but is not
embedded in $F_{R(T)}$. Therefore, the group $F_{R(P)}$ can not be replaced by the terminal
groups in the diagram that the authors constructed for the continuation of the procedure. This is a crucial mistake.

Below is a counterexample. Let $F$ be a free group of rank 4, $F=<a,b,c,d>$. Let $G$ be the double:
$G=<u,v> \, *_{u^3v^3=x^3y^3} \, <x,y>$. $F$ embeds into $G$ by $\nu: F \to G$, where $\nu$ maps the generating
tuple $a,b,c,d$ to $u^4,v^4,x^4,y^4$. However, $G$ has a resolution (fundamental sequence) that maps $G$ onto
a free group of rank 2 (by identifying the images of $<u,v>$ and $<x,y>$), and $\nu(F)$ is mapped onto a (free) subgroup
of rank 2 in that image. In particular, $F$ does not embed into the terminal level of that resolution. In any case,
in contradiction to what the authors claim, 
one can not replace the original limit group $G$, with the embedding $\nu:F \to G$, with the terminal limit group
of the resolution (fundamental sequence) and the image of $F$ in that terminal limit group, which is a
proper quotient of $F$.   

\item"{(27)}" Page 528 line 25: "Therefore, we can further assume that (it) is freely indecomposable  modulo these factors". A
false (critical) reduction. Probably caused by the same reasons that led to the mistake in lines 19-21.

To construct a counterexample let $G$ be the limit group:
$$G \, = \, <u,v> \, *_{[u^3v^3]=t} \, <t,x_1,y_1,x_2,y_2> \, *$$
$$*_{t[x_1,y_1][x_2,y_2]=[z^3,w^3]} \, <z,w>$$
Let $F$ be a free group $<a,b>$, and let $\nu:F \to G$, map the pair $a,b$ to $x_1^5,y_2^7$. $G$ has a resolution
(fundamental sequence) in which $G$ is mapped to the free product $F_2*F_2$, where $<u,v>$ is mapped onto the first factor,
and $<z,w>$ is mapped onto the second factor. $F_2*F_2$ is freely decomposable, but it is certainly not true that to analyze
the Diophantine set that is associated with the embedding $\nu: F \to G$, it is enough to consider the two separate factors
in the free decomposition $F_2*F_2$ (note that we gave an example where $F_2*F_2$ is free, but it's easy to give an example
with freely indecomposable factors).

\item"{(28)}" Page 529  
lines 17-18: "because we add only elements from abelian subgroups". This is completely false. It follows from the 
wrong construction of the induced resolution ("Induced NTQ system") in section 7.12. The counterexample that is
presented in comment
(15) clarifies that in general one must add additional elements, not only elements from "abelian groups" as the
authors claim.




\item"{(29)}" Page 529 line 29: "by levels from top to the bottom". This is wrong. The construction of the 
$induced$ $resolution$ is done iteratively from top to bottom. However, to calculate dimensions, when adding 
QH vertex groups and pegs of abelian vertex groups, and mostly to finally obtain a $firm$ subresolution (definition 4.1
in [Se5]), one has to go from bottom to top.

\item"{(30)}" Page 529 lines -10 to -8: Note that the notion of 
"dimension" of a "fundamental sequence", which is central in the construction of the "tight NTQ envelope",
 has a meaning only for fundamental sequences that satisfy the restrictions that
are listed in sections 7.8 and 7.9 on pages 504-505 (our $well$-$separated$ resolutions). However, these restrictions
are not assumed or used in the construction of "induced NTQ systems" in section 7.12 on page 506, 
and this construction is essential in the construction of "tight NTQ envelope" (i.e., the core resolution). 

\item"{(31)}" Page 529 lines 11 to -8. In these lines, the authors construct their "tight enveloping system" (line -10),
which is the authors' analogue for our $core$ $resolutions$. 
However, the authors slightly modified the construction in section 4 of [Se5], and
what they are actually describing is completely wrong, and can not serve
for the purposes that will be needed in the sequel, i.e., for the purposes that the $core$ $resolutions$ in section
4 of [Se5] was constructed for. 

Needless to say, this mistaken construction immediately leads to false statements in the sequel, 
that can be fixed only if the
construction that is described by the authors is replaced by our actual $core$ $resolutions$.

\item"{(32)}" Page 529 lines -3 to page 530 line 1: "If the dimensions are the same, we can always reorganize the 
levels so that... 
  If all the parameters .. are the same, then $TEnv(S_1)$ has
one level the same as $S_1$". 
This statement of the authors  is valid for our $core$
$resolution$ (theorem 4.13 in [Se5]), and  is indeed a fundamental property of the $core$ $resolution$
and one of its basic properties, but it is completely false for the construction of the "tight enveloping system"
that the authors presented in this page (lines 11 to -8).

Consider the following counterexample (that is based on the counterexample from 2001 [Se8]). Let $F_{2g}$ be a free
group of rank $2g$, $g \geq 2$. Let $L$ be the limit group:
$$L \, = \, <x_1,x_2,x_3> \, *_{x_3[x_1^3,x_2^3]=[u_1,u_2] \ldots [u_{2g-1},u_{2g}]} \, <u_1,\ldots,u_{2g}>$$
The limit group $L$ embeds into the NTQ group (i.e., the $completion$ or the tower) $N$:
$$N \, = \, L \,*_{x_3=[v_1,v_2] \ldots [v_{2g-3},v_{2g-2}]} \, <v_1,\ldots, v_{2g-2}>$$
$F_{2g}$ naturally embeds into $N$ by mapping it onto the (punctured surface subgroup) $<u_1,\ldots,u_{2g}>$. 
We denote this image of $F_{2g}$ in $N$ by $Q$. By the construction of the induced  NTQ
system (section 7.12 of the authors), the induced NTQ system that $Q$ inherits from the NTQ $N$ is:
$$M \, = \,
<u_1,\ldots,u_{2g}> \,
*_{[u_1,u_2] \ldots [u_{2g-1},u_{2g}]=x_3[x_1^3,x_2^3]} \,
<x_1^3,x_2^3,x_3> \,*$$
$$*_{x_3=[v_1,v_2] \ldots [v_{2g-3},v_{2g-2}]} \, <v_1,\ldots, v_{2g-2}>$$
And this "induced NTQ system" $M$  is also the "tight enveloping system" of the authors in this case.
Now, the "dimension" (in the terminology of the authors) of $F_{2g}$ is $2g$, which is exactly the dimension of 
the induced NTQ system $M$. However, the "size" of $M$ is bigger than the size of $F_{2g}$, a contradiction to what
the authors claim on page 529 line -1 (and this is a fatal mistake, as the inequality claimed by the authors is
a fundamental property that is later needed for proving that the procedure for quantifier elimination actually
terminates).

\item"{(33)}" Page 530 lines 1 to 4: "Notice that the dimension of the tight enveloping NTQ fundamental sequence... 
is the same as the maximal  dimension of the corresponding subgroup in the terminal group in the 
enveloping fundamental sequence modulo..." This is completely false.

It is easy to cook a counterexample to this statement using a modification of the counterexample that we used in
comment (32); hence, we don't see a point in writing such an example in detail. This is again a fatal mistake, as
it is an essential property of our $core$ $resolution$ (but not of the "tight enveloping system" 
that the authors
describe), that the procedure for quantifier elimination is based on.

\item"{(34)}" One of the main properties of  our $core$ $resolution$ (in comparison with induced resolution) is that
it is a $firm$ resolution (definition 4.1 in [Se5]). 
This is crucial for the termination of the sieve procedure (that the 
authors borrow in the sequel), and makes our construction of a core resolution
(section 4 in [Se5]) much more complicated. However, the authors don't seem to mention or care about the $firm$ condition
(or property).

\item"{(35)}" Page 530 lines 6-7: Lemma 25 is crucial. 
In the construction that the authors 
present they take "solutions of $F_{R(M_i)}$ minimal with respect to the group of canonical automorphisms 
corresponding to this splitting" (page 529 lines 7-8), and they don't discuss any analogue of our
$auxiliary$ $resolutions$. This is what we do in case there are no parameters,
and it would suffice for a similar lemma with no parameters (see proposition 4.3 in [Se4]). However, taking
these minimal solutions does not suffice in the presence of parameters. Hence, the lemma as stated is false.

\item"{(36)}" Page 530 lines 16-20. The authors are explaining how to find the fundamental sequences that are discussed in
lemma 25 effectively. However, the "algorithm" they give is completely false. 
One can not analyze all the homomorphisms of a tower  (instead of only minimal ones), and then continue the
fundamental sequence after a Diophantine condition is added. This is  misunderstanding of some basic concepts and
a crucial mistake. 

The following is a counterexample - it uses the example in comment (15).
Let $L$ be the limit group:
$$L \, = \, <a_1,a_2> \, *_{[a_1,a_2]=[b_1,b_2]} \, <b_1,b_2> \, *_{[b_1,b_2]=[c_1,c_2]} \, <c_1,c_2>$$
$L$ admits the resolution ("fundamental sequence"
in the terminology of the authors): $L \to \pi_1(S_2)$, where $S_2$ is the fundamental group of an
orientable  surface of genus 2. Hence, $L$ naturally embeds into the $NTQ$ group $N$, that is the completion
of the given resolution: 
$$N \, = \, L \, *_{[b_1,b_2]=[d_1,d_2]} \, <d_1,d_2> \, *_{[d_1,d_2]=[e_1,e_2]} \, <e_1,e_2>$$

Now, we look at a quotient of $N$ by imposing the new relations:
$$d_1(a_1)^{-1},d_2(a_2)^{-1},e_1(b_1)^{-1},e_2(b_2)^{-1}.$$ 
The group that is obtained from the NTQ limit group $N$ by adding these relations has several (maximal) limit quotients.
One of these maximal limit quotients is $L$ itself, where the map from the NTQ group $N$ onto $L$ is the natural
retraction. In our terminology this retraction $L$ is $solid$ with respect to the group that is generated
by $\pi_1(S_2)$, hence, its resolution (fundamental sequence) with respect to the base level of $N$ has a single step
(this quotient of $N$). But $L$ is embedded in this quotient, where the map from $L$ to any terminal limit group
of a graded (relative) resolution (fundamental sequence) that is obtained from homomorphisms in $shortest$ $form$
(minimal in the authors terminology) can not be a monomorphism by proposition 4.3 in [Se4], a contradiction to the
claim of the authors that their "algorithm" produces fundamental sequences (or terminal limit groups of these
sequences) that are obtained from minimal homomorphisms.

\item"{(37)}" Page 530 line 25: "Amalgamate ... along $\bar H$". In general, $\bar H$ is not embedded in a fundamental
sequence that is associated with it. Hence, in such a case, one needs to replace the fundamental sequence that appears 
in top, or take appropriate limit quotients (see comment (16) for a counterexample).

\item"{(38)}" Page 530 lines 27-28: "Consider a fundamental sequence obtained by taking a fundamental sequence 
for ... and pasting to
it a fundamental sequence corresponding to $\bar H$".  
The mistaken construction of the 
"tight enveloping system" that caused severe problems in the construction of the "block NTQ" (i.e., the $anvil$),
causes similar difficulties here. Hence, the construction that the authors describe in
lines 27-28 will work in the sequel, only if their "tight enveloping system" will be fixed to coincide 
with (or to have the same fundamental properties as)
our $core$ $resolution$.

In particular, the counterexamples in comments (15), (32) and (33) show that the construction that is described in Case 1
can not be executed. On line 25 (page 530) it is claimed that the original group is replaced by a proper
quotient, but the counterexample in comment (15) shows that this is not always the case. The counterexample in comment (32) shows that not always there is a reduction in the complexity of the quotient resolution, and the complexity may actually
increase. The counterexample in comment (33) (that we didn't specify in detail) show that Case 2 (lines -11 to -7 on
page 530) can not be executed - this is because the complexity may increase, and if there is an equality in
complexities, this does not mean that the modular groups that are associated with the "tight enveloping system" can
be identified with the original modular group that are associated with what the authors denote by $S_1$ (see line -11).

\item"{(39)}" Page 
531 lines 9-10: "which do not have  sufficient splittings modulo $F_{R(U)}$ and which are terminal groups 
for fundamental sequences modulo $F_{R(U)}$". This is a serious mistake. A fundamental sequence, as the authors suggest, 
may change the values of the original group, hence, it is absolutely forbidden to use such sequences as a preliminary 
step. The authors made a similar mistake in their first step.

\item"{(40)}" Page 531 lines 20-21: "We construct .. modulo the rigid subgroups in the decomposition of the top 
level..". Without using
our $auxiliary$ $resolutions$ (as they appear in definitions 8 and 9 in [Se6]), what the authors do
 does not make any sense and immediately leads to critical mistakes in the sequel.

\item"{(41)}" Page 531 lines -16 to -15: "then using Lemma 11 we can make its size to be smaller than $size(TEnv(W_{n-1}))$". 
Because the construction of the "tight enveloping system" on page
529 is wrong (i.e., doesn't have the fundamental properties of our $core$ $resolutions$), theorem
4.13 in [Se5] is not valid here, and the statement of the authors is completely false. This is indeed a fatal mistake.
See the counterexample in comment (32) that applies here as well.

\item"{(42)}" Page 531 lines -15 to -13: "Similarly to Lemma 25... that the image of H in $E_{n-1}$ ... 
is a proper quotient of 
$E_{n-1}$ or ... does not have a sufficient splitting modulo J". 
Since the authors don't use
our $auxiliary$ $resolutions$ and they consider (multi-graded) resolutions (fundamental sequences) "modulo rigid
subgroups in the top level of $F_{R(L_1^{(n-1)})}$" their conclusion (which is similar to the conclusion of our
proposition 3 in [Se6]) is false. Again, a critical mistake. 

The counterexample in comment (15) can be easily 
modified to give a counterexample to what the authors claim in these lines (i.e., that "the image $H$ ... in the
terminal group ... is a proper quotient" - lines -5 to -14 on page 531).

\item"{(43)}" Page 531 
lines -4 to -3: "We can consider in this case only fundamental sequences that either do not have a maximal
dimension or.. do not have maximal size". The authors do not provide any argument for that. Indeed, because of the wrong
construction of their "tight enveloping system" the statement is false. The statement is correct if the 
"tight enveloping system" is fixed to be our $core$ $resolutions$ (or a replacement with the same 
fundamental properties), and then the  argument can be found in
theorem 4.18 in [Se5].

The counterexample in comment (32) can be easily modified to provide a counterexample to what the authors claim in these 
lines (note that the counterexample in comment (32) 
show that the complexity (or the "size") may increase, and if it remains equal there is no
direct connection between the modular groups of the new "tight enveloping system" and the original modular
groups).

\item"{(44)}" Page 532 
line 10: "modulo the variables of the next level of $L_1^{(n-1)}$".  
Without any analogue of our $auxiliary$ $resolution$ (definitions 8 and 9
in [Se6]) this
doesn't make any sense. Once again, the counterexample in comment (15) demonstrates that the construction that is
described in these lines can not work.

\item"{(45)}" Page 532 lines -10 to -1: 
Because of the wrong construction of "tight enveloping systems" on page 529, theorem 38
as stated is wrong. 

Each of the counterexamples in comments  (32) and (33) gives a counterexample to the
conclusion of theorem 38. In fact, it is not difficult to use these counterexamples to construct a procedure
of the type that is described by the authors that never terminate (somewhat similar to the non-terminating
procedure that we presented in 2001 [Se8]).
 
The theorem 
can be made correct only if the "tight enveloping system" is replaced by our
$core$ $resolutions$ or an object with the same fundamental properties. Even after such a replacement, 
the arguments and the constructions that the authors use for the proof of theorem 38 need to be (significantly)
modified to agree with the constructions that are presented in the general step of the sieve procedure in
[Se6].

\item"{(46)}" Page 533 lines 9 to 16. The argument regarding the size of the "tight enveloping system" that the authors
present is completely false (because of the wrong construction of the authors "tight enveloping system" on page 529). 
It would work for our $core$ $resolution$, according to theorem 4.13 in [Se5]. For a counterexample, see the 
counterexample in comment (32).

\item"{(47)}" Page 
541 line 18: "generic family of solutions". The authors have defined "generic families" in definition 23, page 508
lines 15-25. 
However, this definition of "generic family" can be used only in the procedure for the analysis of sentences with
2 quantifiers (the authors' construction of the tree $T_{AE}(G)$ in section 9 that starts on page 514 line 8). In
the sieve procedure, i.e., the procedure that the authors borrowed for the analysis of formulas with more than 2
quantifiers (the one that is used in their "projective images" section), the notion of a "generic family" as 
 defined in definition 23, is false. For that procedure, one needs to consider reduced modular groups, and technically
to use any form of Merzlyakov theorem in the current context the authors should have borrowed our notion of
$framed$ $resolutions$ (definition 5 in [Se6]). Without framed resolutions, the use of "generic families" and 
general forms of Merzlyakov's theorem is completely wrong in the current context.

As an example let $G=\pi_1(S_2)$, the fundamental group of a surface of genus 2. Let $\nu: G \to Z_2$ (where $Z_2$
is the group with two elements), and let $H$ be the kernel of this map. Then there exists a fixed element in $h \in H$,
so that for a "generic" homomorphism $f:H \to F_k$, $f(h)$ does not have a square root. On the other hand for every
homomorphism $v:G \to F_k$, $v(h)$ does have a square root. Hence, (for first order considerations) 
the "genericity" notions are different, and the one
 that is used by the authors in the sequel is the wrong one.  

\item"{(48)}" Page 542 line 16: "the sequence of proper projective images stabilizes". 
The authors
apply their
theorem 38 on page 532  to get these stable "projective images". However, as we have already commented in (45),
the conclusion of theorem 38, which is supposed to be one of the authors' main results, is false (this is mainly because
of the wrong construction of their "tight enveloping system", see comments (32) and (33)).

\item"{(49)}" Because of the authors' wrong perception of our $core$ $resolution$ (their
"tight enveloping system"), and their wrong interpretation  of "generic families" in the context of the procedure
for quantifier elimination (our sieve
procedure), the authors  have no way (i.e., a notion or a construction) to distinguish between our
$core$ $resolution$ and our $penetrated$ $core$ $resolution$ (definition 4.20 in [Se4]). Both of these resolutions,
the $core$ and the $penetrated$ $core$ are essential in obtaining any form of stability of what the authors
call "block NTQ" (our $anvils$) and their associated fundamental sequences (our $developing$ $resolutions$). Without
these notions (or rather constructions), and without examining their appropriate $induced$ $resolutions$
(see definition 9 in [Se6]), there is no real meaning to the sentence "projective images stabilize". To see how
the sentence of the authors could be made precise (with all the notions and procedures at hand) see the statement 
and the proof of proposition 23 in [Se6].

\item"{(50)}" Page 541 lines 20-24: "fundamental sequences induced by the subgroup of the enveloping 
group generated by... We call this group a second principal group and consider projective images for these fundamental
sequences". The authors are trying 
to define our $sculpted$ $resolution$ as it appears in part (4) of the first step of the sieve procedure,
and parts (6) and (7) of the general step of the sieve procedure in [Se6]. However, what they do write
is completely
wrong 
(i.e., one can not work or continue the procedure with the construction that the authors suggest, and hope for
termination). 

The authors perform things in reverse order - once one fixes the (second) $algebraic$ $envelope$, one should
first look at its sequence of $core$ $resolutions$ ("tight enveloping system") and then construct 
a resolution that is composed of a sequence of induced resolutions.

\item"{(51)}" Page 541 lines 23-24: "second principal group" is supposed to be the authors' analogue of our 
$second$ $sculpted$
$resolution$, as it appears in part (4) of the first step of the sieve procedure. Probably because of the
mistaken construction of their "tight enveloping system" on page 529, the authors
do not realize that beside the $sculpted$ $resolution$, one needs to construct, keep, and examine the structure
of a related resolution, the $penetrated$ $sculpted$ $resolution$ (see parts (6)-(7) of the general step of the sieve
procedure in [Se6]). Without $penetrated$ $core$ $resolutions$ and $penetrated$ $sculpted$ $resolutions$, the
procedure that the authors describe has no reason to terminate.

\item"{(52)}" Page 541 lines -11 to -9: 
"When the chain of projective images of these sequences of increasing depth (and fixed thickness) stabilizes, 
we add new variables, increase thickness, and so on". This is wrong.

In the general step of the sieve procedure (the procedure for quantifier elimination)
in [Se6], it is impossible (conceptually) to "know" when the constructions
that are associated with a certain algebraic envelope stabilize. What we do 
in the general step (parts (6)-(7)) is 
check the objects that are associated with every algebraic envelope (sculpted, penetrated sculpted, and developing
resolutions, core and penetrated core resolutions, Carriers) at every step, going from the first algebraic
envelope to the last (the one with maximal $width$), and for the first one in which there is a change we
modify the constructions that are associated with it, and basically forget or remove the higher width algebraic
envelopes. If there is no change in any of them we add an algebraic envelope of a bigger width.

A major principle of the sieve procedure is that not only  the constructions that are associated with 
algebraic envelopes are changed along the procedure, but also the algebraic envelopes themselves are
changed! only when we examine an infinite path of the procedure we can talk about stable algebraic envelope with
its stable associated objects (see proposition 23 and definition 24 in [Se6]).

Therefore, what the authors write in lines -14 to -9 on page 541 does not agree with the basic concepts of the
sieve procedure in [Se6]. From our point of view, it's not only a mistake, but a basic misunderstanding that doesn't 
allow one to construct our sieve  procedure, or in fact any procedure that is supposed to analyze Diophantine or
more generally definable sets.

\item"{(53)}" Page 542 line 6: "minimal values from some NTQ system". What does minimal have to do in this context?

\item"{(54)}" Page 543 lines 1 to 3: "As we did in the proof of theorem 11, we construct for $K_E$... such that 
$K_E$ is not conjugated into a fundamental group of a proper subgraph of the JSJ decomposition D ...". This 
statement is wrong, but it is supposed
to be the authors' translation of theorem 32 in [Se6]. 

\item"{(55)}" Page 543 lines 8-9: "Therefore the quadratic system ... has smaller size than the quadratic system of..". This is 
precisely proposition 37 in [Se6]. The proof of this proposition is not as simple as the authors claim in lines
5-8, and requires  significant modifications to the argument that was used to prove theorems 2.5 and 2.9 in [Se3]
(i.e., what the authors call theorem 11).

\item"{(56)}" Page 544 lines 1-14. 
The exact statement of lemma 29 is still wrong (it is identical to the statement of theorem 1.33 in [Se5] that
is false in this context). The  authors use the term "generic family" in part (1) of the lemma
(line 6). They define a "generic family" in definition 23 on page 508 lines 15-25. However, this definition of
"generic family" (i.e., our $graded$ $test$ $sequence$) is invalid here. One has to define
$framed$ $resolutions$ (as they appear in definition 5 in [Se6]), and prove a form of a generalized
Merzlyakov theorem, in order to make the notion "generic family" precise 
in the context of the lemma (see the difference between the formulations of lemma 6 in [Se6] and theorem
1.33 in [Se5]). For an example for the difference between the different notions of "generic families" see comment (47).

\item"{(57)}" Page 546 line -10 to -9. Theorem 41  is identical with
theorem 7 in [Se7]. 
The formulation of theorem 7 in [Se7] has a mistake, that was found and fixed by C. Perin [Pe]. Exactly the same
mistake exists in the formulation of theorem 41 of the authors. 

\endroster


\vskip 2em

\end